\documentclass[12pt,reqno]{amsart}
\usepackage[top=2cm,bottom=2cm,left=1.5cm,right=1.5cm]{geometry}
\usepackage[T1]{fontenc}
\usepackage{lmodern}
\usepackage{amsmath,amssymb,amsthm,mathtools,bm}
\usepackage{graphicx,xcolor,enumitem,microtype}
\usepackage[colorlinks=true,citecolor=blue,linkcolor=blue,urlcolor=blue]{hyperref}
\allowdisplaybreaks
\newif\ifshowchanges
\showchangestrue
\newcommand{\rev}[1]{{\ifshowchanges\color{blue}\fi#1}}

\newtheorem{theorem}{Theorem}[section]
\newtheorem{proposition}[theorem]{Proposition}
\newtheorem{lemma}[theorem]{Lemma}
\newtheorem{corollary}[theorem]{Corollary}
\theoremstyle{definition}
\newtheorem{definition}[theorem]{Definition}

\theoremstyle{remark}
\newtheorem{remark}[theorem]{Remark}

\newcommand{\N}{\mathbb{N}}
\newcommand{\Nzero}{\mathbb{N}_0}
\newcommand{\R}{\mathbb{R}}

\newcommand{\E}{\mathbb{E}}
\newcommand{\Pp}{\mathbb{P}}
\newcommand{\e}{\vec{\bm e}}
\newcommand{\x}{\vec{\bm x}}
\newcommand{\z}{\vec{\bm z}}
\newcommand{\kvec}{\vec{\bm k}}
\newcommand{\bvec}{\vec{\bm b}}

\newcommand{\hvec}{\vec{\bm h}}
\newcommand{\ivec}{\vec{\bm i}}
\newcommand{\uvec}{\vec{\bm u}}

\newcommand{\zerovec}{\vec{\bm 0}}
\newcommand{\Mvec}{\vec{\bm M}}
\newcommand{\Uvec}{\vec{\bm U}}
\newcommand{\X}{\vec{\bm X}}

\newcommand{\V}{\mathcal{V}}
\newcommand{\Kc}{\mathcal{K}}

\newcommand{\diag}{\operatorname{diag}}

\newcommand{\dd}{\,\mathrm{d}}
\newcommand{\ind}{\mathbf{1}}
\newcommand{\abs}[1]{\left\lvert #1\right\rvert}
\newcommand{\norm}[1]{\left\lVert #1\right\rVert}
\newcommand{\binomulti}[2]{\binom{#1}{#2}}
\newcommand{\eqdist}{\mathrel{\stackrel{d}{=}}}
\DeclareMathOperator{\Mult}{Mult}

\title[Multivariate Bernstein duality]{Multivariate Bernstein duality for tournament Wright--Fisher diffusions}

\author[O. Angtuncio Hern\'andez]{Osvaldo Angtuncio Hern\'andez}
\address{Universidad de las Am\'ericas Puebla, Departamento de Actuar\'ia, F\'isica y Matem\'aticas, Puebla, M\'exico}
\email{osvaldo.angtuncio@udlap.mx}


\hypersetup{
pdftitle={Multivariate Bernstein duality for tournament Wright-Fisher diffusions},
pdfauthor={Osvaldo Angtuncio Hernandez}
}

\newcommand{\reviewfigure}[2]{
\IfFileExists{#2}{
\includegraphics[width=#1]{#2}
}{
\fbox{\begin{minipage}[c][3cm][c]{.85\linewidth}
\centering
\textcolor{blue}{Corrected figure not embedded in this text-only copy.}
\par\medskip
\small\texttt{\detokenize{#2}}
\end{minipage}}
}
}

\begin{document}

\begin{abstract}
We give a probabilistic proof of a multivariate Bernstein duality for a multivariate Wright--Fisher frequency process with selection encoded by a finite tournament.
The construction includes nontransitive interactions, such as rock--paper--scissors selection.
Its dual is a pure-jump process of coefficient functions in the multinomial Bernstein basis.
Selection increases the Bernstein degree through a winner--loser averaging operator, whereas neutral resampling decreases it through a duplication operator.
For the proof, we construct a continuous forward process independently of duality, use duality to establish uniqueness, and obtain a weak stochastic differential equation representation and invariance of the simplex faces.
This provides a sampling-based treatment of the tournament specialization of multitype Bernstein duality, extending an open problem in the binary Kingman setting of Cordero, Hummel, and Schertzer (2022).
\end{abstract}

\maketitle
\enlargethispage{3pt}

\section{Introduction}\label{sec:introduction}

In the neutral Wright--Fisher model, a population of fixed size $N$ reproduces in discrete generations, and each individual carries one of two possible types.
Each offspring independently chooses a parent uniformly from the previous generation and inherits its type.
The model can be studied forward in time through type frequencies or backward in time through the ancestry of a sample.
After rescaling time, the frequency process converges to the Wright--Fisher diffusion; see \cite[Chapter~10]{EthierKurtz1986}.
The corresponding ancestral limit is Kingman's coalescent; see \cite{MR671034,MR1880231}.

Duality connects these forward and backward descriptions.
In the limit, the frequency process solves the equation
\[
    \mathrm dX_t=\sqrt{X_t(1-X_t)}\,\mathrm dB_t,
    \qquad X_0=x\in[0,1],
\]
where $B$ is a standard Brownian motion.
On the other hand, the neutral dual is the line-counting process $L$ of Kingman's coalescent, with transition $m\mapsto m-1$ at rate $\binom m2$.
The classical moment identity between them is
\[
    \E_x[X_t^n]=\E_n[x^{L_t}],
    \qquad t\geq0,\quad n\in\N.
\]
Under genic selection, the backward object is the ancestral selection graph (ASG) of Krone and Neuhauser \cite{NeuhauserKrone1997,KroneNeuhauser1997}.
General discussions of duality and population-genetic models may be found in \rev{\cite{Etheridge2011,JansenKurt2014}}.

Cordero, Hummel, and Schertzer \cite{CorderoHummelSchertzer2022} introduced Bernstein duality for two-type $\Lambda$-Wright--Fisher processes with polynomial frequency-dependent selection.
These processes arise as frequency limits of population models with selective and neutral reproduction events.
At a selective event, an individual is assigned a random collection of potential parents, and a prescribed selection rule, depending on their types, determines which parent replaces the individual.
At a neutral reproduction event, an individual produces offspring that replace a fraction of the population.
Their Theorem~2.14 establishes a \emph{Bernstein duality} between the limiting frequency process and a coefficient-valued Markov process driven by branching and coalescence.
The latter is described in Definitions~2.11--2.12, while Proposition~2.13 identifies its coefficients with those of the ancestral selection polynomial (ASP).
Informally, the ASP evaluated at $x\in[0,1]$ and backward time $t$ is the conditional probability that all sampled roots have a prescribed type, given the ancestral line-count history, when the potential ancestors at backward time $t$ receive independent types with frequency $x$.
Finally, Question~(9) in Section~2.10 explicitly asks whether this duality can be extended to higher-dimensional frequency processes, motivated by nontransitive interactions between types.
In this paper we show that such an extension is possible, under a binary selection and reproduction mechanism.

The authors prove in Equation~(1.4) the Bernstein duality, which is
\begin{equation}\label{eq:CHS-intro-duality}
    \E_x[X_t^n]
    =\E_{e_{n+1}}\left[
        \sum_{j=0}^{L_t}V_t(j)\binom{L_t}{j}
        x^j(1-x)^{L_t-j}
    \right],
\end{equation}
where $V_0=e_{n+1}:=(0,\ldots,0,1)$, $L_t:=\dim(V_t)-1$, $n\in\N$ and $x\in[0,1]$.
The coefficients $(V_t(j))_{j\leq L_t}$ in the Bernstein basis of degree $L_t$ average over ancestral graphs with the same line-count history.
We extend this binomial representation to multinomial coefficients and arbitrary finite tournaments, while restricting reproduction to the binary case.

Gonz\'alez Casanova and Smadi \cite{GonzalezCasanovaSmadi2020} construct a more general multitype forward model with colouring rules, mutation, Wright--Fisher noise, and extreme reproductive events.
Their Section~4.3, Equation~(4.4), includes pairwise competitions with prescribed winning probabilities, and their Section~4.4, Equation~(4.5), gives rock--paper--scissors selection.
Their Remark~4.2 explains why the classical moment-duality method used for transitive rules does not apply to general nontransitive colouring rules.
The present paper supplies an explicit coefficient-valued dual for the deterministic tournament and Brownian-resampling subclass of that framework.
It does not generalize their full forward construction.

Several later works extend the one-dimensional theory in other directions.
Cordero and V\'echambre \cite{CorderoVechambre2023} study selection and mutation in a one-sided random environment.
Cordero, Hummel, and V\'echambre \cite{CorderoHummelVechambre2025} treat fixation and coexistence under opposing environmental effects, and \cite{CorderoHummelVechambre2026} uses Bernstein duality to study long-term behaviour with polynomial selection, mutation, and environmental effects.
Cordero and V\'echambre \cite{CorderoVechambre2025Distance} use Bernstein and Siegmund duals for quantitative ergodic estimates and open-set recurrence.
Our extension instead replaces the scalar frequency by a simplex-valued frequency vector and the binomial basis by the multinomial basis.
We retain binary deterministic selection and Kingman resampling; mutation, higher-order selective interactions, and non-Brownian $\Lambda$-reproduction events are not included.

Very recently, Hummel \cite[Theorem~2.3]{Hummel2026Multitype} has established Bernstein duality for a class of multitype $\Lambda$-Wright--Fisher processes with finite polynomial colouring mechanisms and mutation.
The author's main theorem contains the present tournament diffusion as the specialization $\Lambda=\gamma\delta_0$, with deterministic binary selection and no mutation.
The contribution of the present work is therefore complementary rather than an extension in generality.
Whereas Hummel proves the generator identity by several multinomial coefficient calculations, we derive it through explicit probabilistic sampling identities: independent colouring reduces selective events to a winner--loser experiment, while neutral resampling is represented by duplication of a uniformly sampled type.
These constructions give a direct interpretation of the coefficient operators in terms of sampled individuals.
Thus, the value of the present treatment lies in its transparent probabilistic mechanism and its explicit connection between tournament competitions
and Bernstein coefficients.
Another new contribution of the paper is that we also describe the forward frequency dynamics as the evolution of a population, and as a stochastic differential equation. This helps to explain how selective competitions and neutral reproductions contribute to their generator.

In our tournament model, at each selective event, the target
receives the winning type of the challenger and the target.
If the target wins, neither individual's type changes.
Moreover, we add neutral reproductions in which an individual reproduces and replaces another individual's type by a copy of its own, irrespective of the tournament.
A biological motivation for tournament selection is the classical three-morph mating system of the side-blotched lizard \emph{Uta stansburiana}; see \cite{SinervoLively1996,VeraRamirezEtAl2012}.
Male throat colour is associated with alternative reproductive strategies: orange males defend large territories and have an advantage over blue males, blue males guard smaller territories effectively against yellow males, and yellow males exploit opportunities for sneaking matings within the territories of orange males.
The resulting orange--blue--yellow--orange cycle motivates a nontransitive tournament rather than a fixed ranking of types.
Our model idealizes these relative reproductive advantages by deterministic pairwise competition, while neutral reproduction introduces demographic randomness.
The types represent reproductive strategies within one species, not different species.

The paper is organized as follows. Section~\ref{sec:main-result} defines the forward process and the two ancestral coefficient operators and states our main result, Theorem~\ref{thm:bernstein-duality}.
Section~\ref{sec:bernstein} establishes the sampling and generator identities.
In particular, Lemma~\ref{lem:winner-loser-identity} shows that deleting the loser of two uniformly marked individuals in an independent sample has the same joint law as adjoining the winner of two independent types to an independent smaller sample.
Lemma~\ref{lem:duplication-identity} gives the corresponding identity for duplication.
Subsection~\ref{sec:coefficient} proves nonexplosion and finite-time degree bounds.
Section~\ref{sec:duality} combines these ingredients to prove duality and then deduces uniqueness of the forward martingale problem.
Appendix~\ref{app:forward} independently constructs a continuous forward solution and establishes its stochastic differential equation representation and some other properties.

One of the main applications of our duality is that it extends \eqref{eq:CHS-intro-duality} to every mixed moment as
\begin{equation}\label{eq:moment-dualityIntro}
    \E_{\x}\left[\prod_{j=1}^dX_j(t)^{k_j}\right]
    =\binomulti n{\kvec}^{-1}
    \E_{\delta_{\kvec}}\left[
       \sum_{\hvec\in\Kc_{L_t}}
       V_t(\hvec)\binomulti{L_t}{\hvec}\x^{\hvec}
    \right].
\end{equation}
Here $n\geq1$, $\kvec\in\Nzero^d$ satisfies $\sum_{j=1}^d k_j=n$, $\delta_{\kvec}$ is the indicator coefficient function of this count vector, and $L_t$ is the current number of potential ancestral lines.
The random array $(V_t(\hvec))_{\hvec\in\Kc_{L_t}}$
gives the coefficients of a polynomial in the multinomial
Bernstein basis of degree $L_t$.
These objects are defined in Subsection~\ref{subsec:main-coefficients}.

In the recent expository article \cite{AngtuncioEslavaFittipaldi2026}, together with Laura Eslava and Mar\'ia Clara Fittipaldi, we develop two  finite-population scaling limits of the population with tournament selection.

\section{Model and main result}\label{sec:main-result}\label{sec:forward}

We begin by describing the population dynamics from both the forward perspective of type frequencies and the backward perspective of ancestry.
Although our focus is on the infinite-population model, we first describe the finite-population replacement mechanisms, which motivate the limiting processes.

\subsection{Forward evolution: tournament frequencies}\label{subsectionForwardEvolution}

Write $\N=\{1,2,\ldots\}$, $\Nzero=\N\cup\{0\}$, and $[m]=\{1,\ldots,m\}$, with $[0]=\varnothing$.
Fix $d\geq2$, let $\e_1,\ldots,\e_d$ be the canonical basis of $\R^d$, and set
\[
    \Delta_d:=\left\{\x\in[0,1]^d:\sum_{j=1}^d x_j=1\right\}.
\]
A tournament matrix $T=(T_{ij})_{i,j\in[d]}$ satisfies $T_{ii}=0$, $T_{ij}\in\{0,1\}$, and $T_{ij}+T_{ji}=1$ for $i\neq j$.
The convention $T_{ij}=1$ means that type $i$ defeats type $j$.
Write $A:=T-T^{\mathsf T}$; thus $A$ is skew-symmetric and $A_{ij}\in\{-1,1\}$ off the diagonal.

\emph{Competitions and reproductions in a finite population.}
Fix a population size $N\in\N$ with $N\geq2$.
There are $N$ labelled individuals, with types $i_u\in[d]$, $u\in[N]$.
They evolve in continuous time through competitions and neutral reproductions.
For each ordered pair of distinct individual labels $(U_1,U_2)$, use independent Poisson clocks with respective rates $\lambda_N,\gamma_N\geq0$, omitting any clock whose rate is zero.
All clocks are independent of the initial types and of one another.
At a \emph{competition}, the target $U_2$ receives the winning type of $i_{U_1}$ and $i_{U_2}$; if the types agree, that common type is retained.
At a \emph{neutral reproduction}, the target $U_2$ copies the type of the source $U_1$.
Only the type at $U_2$ can change, so population size remains $N$.
This is the binary tournament version of the individual construction in \cite[Section~2.1]{CorderoHummelSchertzer2022}.

The frequency vector $\X^N$ has coordinates
$X_j^N=N^{-1}\sum_{u=1}^N\ind_{\{i_u=j\}}$
and state space $\Delta_d\cap N^{-1}\Nzero^d$.
Every transition that changes the frequency vector has the form
\begin{equation}\label{eq:intro-finite-rates}
    \x\longmapsto\x+\frac{\e_i-\e_j}{N},
    \qquad i\neq j,\quad x_ix_j>0.
\end{equation}
Such a transition occurs when a target of type $j$ is replaced by type $i$, either through neutral reproduction by a source of type $i$, or through a competition in which a type-$i$ challenger defeats the type-$j$ target (so $T_{ij}=1$).
Thus $i$ denotes the new type of the target and $j$ its previous type.

\emph{Competitions and reproductions in an  infinite population.}
Here we define the limiting process through the operator
\begin{equation}\label{eq:forward-generator}
    \mathcal Lf(\x)
    :=\lambda\sum_{j=1}^d x_j(A\x)_j\partial_j f(\x)
     +\frac\gamma2\sum_{i,j=1}^d
           x_i(\delta_{ij}-x_j)\partial_{ij}f(\x).
\end{equation}
We allow $\lambda,\gamma\in[0,\infty)$; the diffusion regime has $\gamma>0$, whereas $\gamma=0$ gives a deterministic frequency flow.
Here $\delta_{ij}$ is the Kronecker delta, and $C^2(\Delta_d)$ denotes restrictions to $\Delta_d$ of functions that are twice continuously differentiable on a neighbourhood of the simplex.
The first term describes the effect of competitions, and the second describes random changes caused by reproduction.
For later use, put $b_j(\x):=\lambda x_j(A\x)_j$ and $\bvec=(b_j)_{j\in[d]}$.
Corollary~\ref{cor:forward-wellposed} proves that the continuous martingale problem for this operator is well posed.

\begin{figure}[tbp]
\centering
\reviewfigure{.83\textwidth}{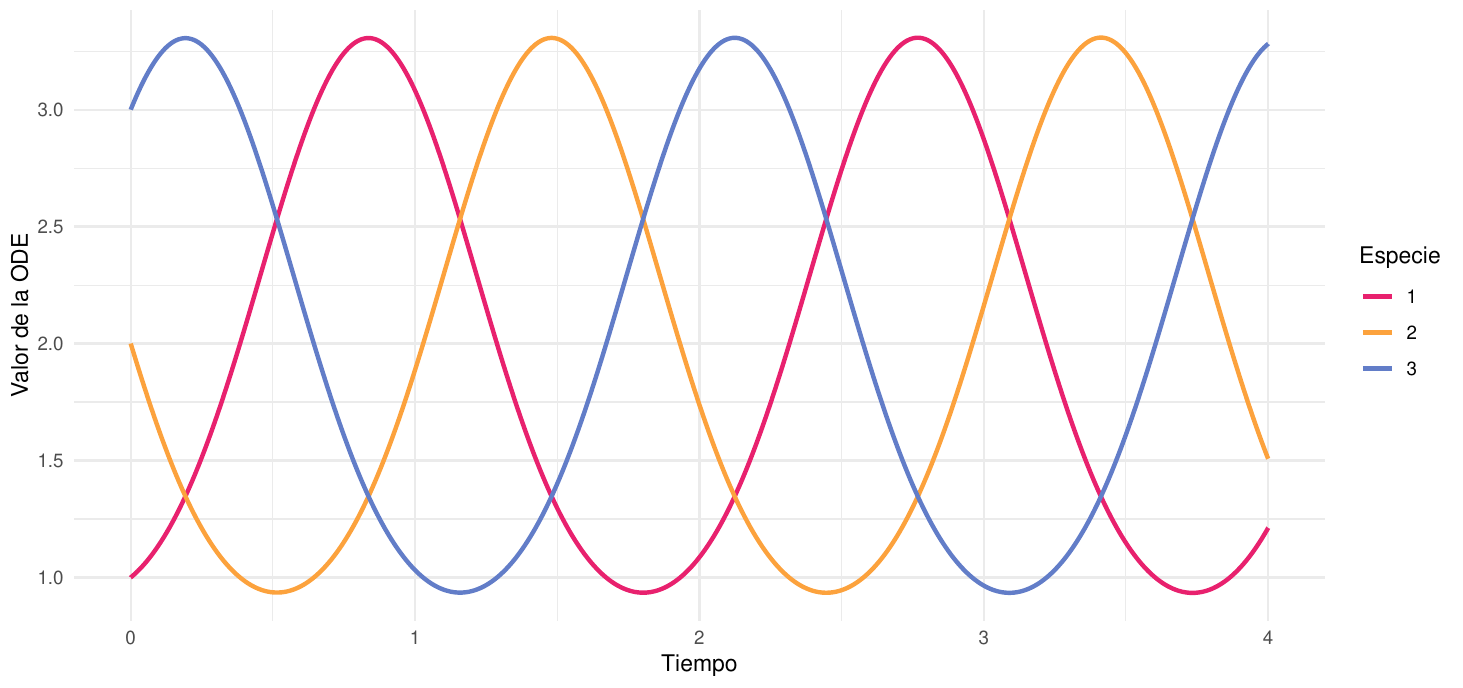}
\caption{Numerical solution of the deterministic cyclic frequency equation with $\x_0=(1,2,3)/6$.
The orange, yellow, and blue curves represent types $1,2,3$, respectively.
The frequencies sum to one; the product $x_1x_2x_3$ is constant along the deterministic trajectory.}
\label{figODE}
\end{figure}

\begin{figure}[tbp]
\centering
\reviewfigure{.83\textwidth}{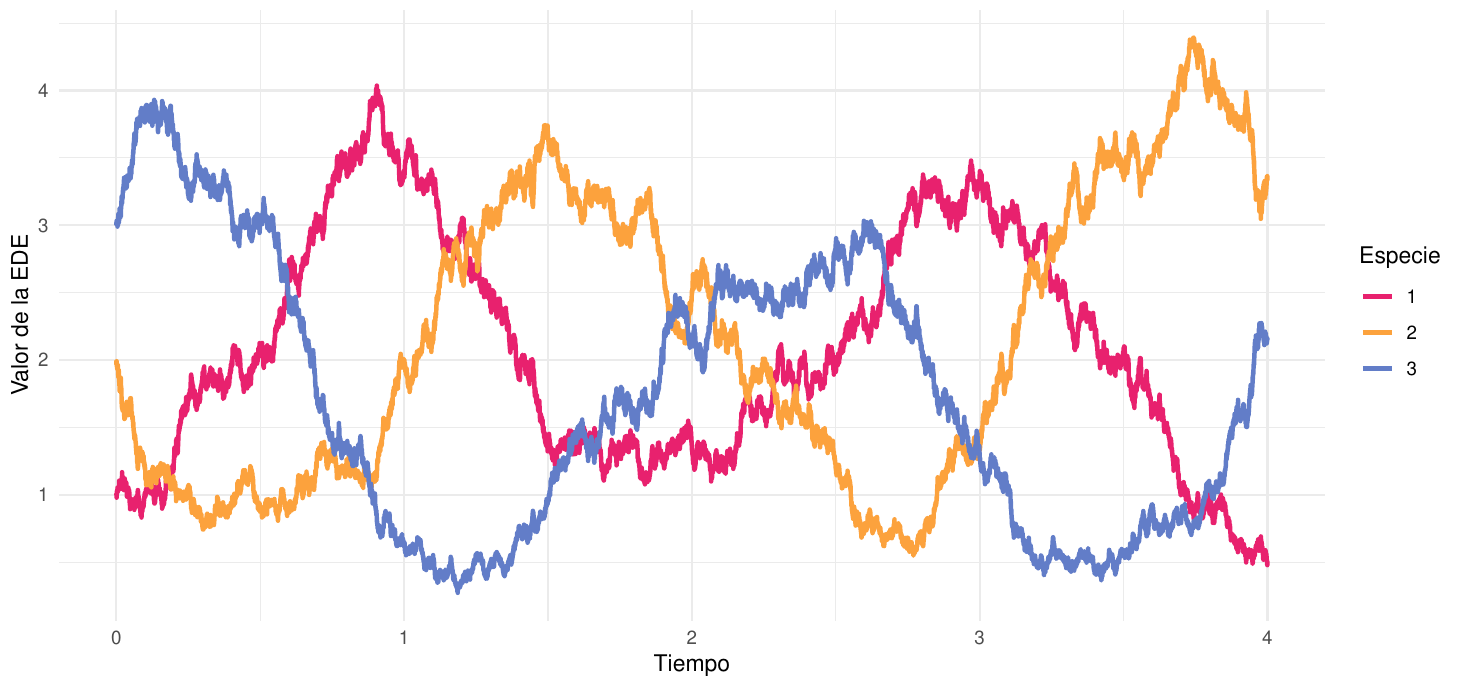}
\caption{One trajectory that approximates the tournament Wright--Fisher diffusion, with $\x_0=(1,2,3)/6$.
The colours and line styles have the same meaning as in Figure~\ref{figODE}.
}
\label{figEDENonFixation}
\end{figure}

For the lizard-inspired example, label the orange, yellow, and blue strategies by $1,2,3$, respectively.
With the convention that $T_{ij}=1$ means that type $i$ defeats type $j$, the cyclic interaction is encoded by
\[
T=\begin{pmatrix}
0&0&1\\
1&0&0\\
0&1&0
\end{pmatrix}.
\]
Thus orange defeats blue, yellow defeats orange, and blue defeats yellow.
A selective replacement represents a reproductive advantage rather than a literal change in the throat colour of an adult.
In the absence of neutral resampling ($\gamma=0$), the frequency process is deterministic and satisfies
\[
\dot x_1=\lambda x_1(x_3-x_2),
\qquad
\dot x_2=\lambda x_2(x_1-x_3),
\qquad
\dot x_3=\lambda x_3(x_2-x_1),
\]where $\dot x_j=\mathrm dx_j/\mathrm dt$. 
Figures~\ref{figODE} and \ref{figEDENonFixation} illustrate this deterministic system and a finite-population approximation to the diffusion, respectively.
In the figures we use a particular rescaling for the state space and time. 

\subsection{Backward ancestry: coefficient functions and ancestral operators}
\label{subsec:main-coefficients}

We now introduce the backward process associated with the limiting frequency process whose generator is given in \eqref{eq:forward-generator}.

\emph{Multitype Bernstein polynomials.}
For $n\in\Nzero$, define
\begin{equation}\label{eq:multi-index-set}
    \Kc_n:=\left\{\kvec\in\Nzero^d:
    \abs{\kvec}:=\sum_{j=1}^d k_j=n\right\}.
\end{equation}
In particular, $\Kc_0=\{\zerovec\}$.
For $\kvec\in\Kc_n$ and $\x\in[0,\infty)^d$, write
\begin{equation}\label{eq:multinomial-notation}
    \binomulti{n}{\kvec}:=\frac{n!}{k_1!\cdots k_d!},
    \qquad
    \x^{\kvec}:=\prod_{j=1}^d x_j^{k_j},
\end{equation}
with $0!=1$ and $0^0=1$.
For $\kvec\in\Kc_n$, define the multinomial Bernstein polynomial
\begin{equation}\label{eq:bernstein-basis}
    B_{\kvec,n}(\x):=\binomulti n{\kvec}\x^{\kvec}.
\end{equation}
Let $\mathcal P_n(\Delta_d)$ be the restrictions to $\Delta_d$ of polynomials of total degree at most $n$.
At a fixed degree, the following lemma says that Bernstein coefficients specify exactly the same information as such a polynomial; see \cite[Section~3, equations~(7)--(8)]{Titi2017-12Matri-40426}.

\begin{lemma}[The multinomial Bernstein basis]\label{lem:bernstein-basis}
For each $n\in\Nzero$, the family $(B_{\kvec,n})_{\kvec\in\Kc_n}$ is a basis of $\mathcal P_n(\Delta_d)$.
Thus every $f\in\mathcal P_n(\Delta_d)$ has a unique coefficient function $v_f^{(n)}:\Kc_n\to\R$ satisfying
\[
    f(\x)=\sum_{\kvec\in\Kc_n}v_f^{(n)}(\kvec)B_{\kvec,n}(\x),
    \qquad \x\in\Delta_d.
\]
The indexed array $(v_f^{(n)}(\kvec))_{\kvec\in\Kc_n}$ is its degree-$n$ multitype Bernstein coefficient vector.
\end{lemma}

A labelled population of size $m\geq1$ has individual labels in $[m]$.
A type assignment, or colouring, is a vector $\ivec=(i_1,\ldots,i_m)\in[d]^m$, where $i_u$ is the type of individual $u$.
Its \emph{degree vector} is
\begin{equation}\label{eq:count-vector-definition}
    \operatorname{deg}(\ivec)
    :=\sum_{u=1}^m\e_{i_u},
    \qquad
    [\operatorname{deg}(\ivec)]_j
    =\sum_{u=1}^m\ind_{\{i_u=j\}}.
\end{equation}
Thus $\operatorname{deg}(\ivec)=\kvec$ means that exactly $k_j$ individuals have type $j$, for every $j\in[d]$.
The empty colouring has degree vector $\zerovec$.
For $\x\in\Delta_d$, saying that a random type $I$ has law $\x$ means $\Pp(I=j)=x_j$ for $j\in[d]$.
Let $I_1,\ldots,I_m$ be independent random types with this law and define
\begin{equation}\label{eq:multinomial-realization}
    \Mvec_m(\x):=\sum_{u=1}^m\e_{I_u},
    \qquad
    \Pp(\Mvec_m(\x)=\kvec)
    =\binomulti m{\kvec}\x^{\kvec},\quad \kvec\in\Kc_m.
\end{equation}
Thus $\Mvec_m(\x)$ has a multinomial law $\Mult(m,\x)$.
For $m=0$, the sum is empty and $\Mvec_0(\x)=\zerovec$ almost surely.
When independence from other variables is specified below, this notation denotes a fresh multinomial vector with the indicated law.
For $n\geq1$, let
\begin{equation}\label{eq:coefficient-spaces}
    \V_n:=\{v:\Kc_n\longrightarrow\R\},
    \qquad
    \V:=\bigsqcup_{n\geq1}\V_n.
\end{equation}
An element $v\in\V_n$ is a coefficient function on $\Kc_n$.
For $v\in\V_n$, set
\begin{equation}\label{eq:coefficient-domain-index}
    d(v):=n,
    \qquad\text{so that }\operatorname{dom}(v)=\Kc_{d(v)}.
\end{equation}
The integer $d(v)$ is its \emph{domain index}, also called its Bernstein degree.
Note that it is distinct from the fixed number $d$ of types.
Because the states are functions, we keep the notation $v$ and $V_t$ without vector accents.
For $v\in\V_n$, define
\begin{equation}\label{eq:coefficient-sup-norm}
    \norm{v}_\infty:=\max_{\kvec\in\Kc_n}\abs{v(\kvec)}.
\end{equation}
The Bernstein evaluation of a coefficient function $v\in\V_n$ is
\begin{equation}\label{eq:H-definition}
    H(v,\x):=\sum_{\kvec\in\Kc_n}
    v(\kvec)B_{\kvec,n}(\x)=\E[v(\Mvec_n(\x))],
\end{equation}
where the expectation follows from \eqref{eq:multinomial-realization}.
The probabilistic coefficient functions form the invariant subset
\begin{equation}\label{eq:probabilistic-coefficients}
    \V^{[0,1]}:=\bigsqcup_{n\geq1}\{v\in\V_n:0\leq v(\kvec)\leq1
        \text{ for all }\kvec\in\Kc_n\}.
\end{equation}

\emph{The selection and coalescent operators.}
We now define the operators corresponding to the two transitions of the backward process.
Simple probabilistic interpretations are proved in Section~\ref{sec:bernstein}.
For two types $i,j\in[d]$, define their tournament winner and loser by
\begin{equation}\label{eq:winner-loser-map}
    w(i,j):=\begin{cases}
        i,&i=j\text{ or }T_{ij}=1,\\
        j,&T_{ji}=1,
    \end{cases}
    \qquad
    \ell(i,j):=\begin{cases}
        i,&i=j\text{ or }T_{ji}=1,\\
        j,&T_{ij}=1.
    \end{cases}
\end{equation}
These are the winner and loser of the competitions described in Section~\ref{sec:forward}.
For $n\geq1$ and $\kvec\in\Kc_{n+1}$, define
\begin{equation}\label{eq:q-winner-loser}
    q_{ij}^{(n+1)}(\kvec):=
    \begin{cases}
        \displaystyle\frac{2k_ik_jT_{ij}}{n(n+1)},&i\neq j,\\[1.2ex]
        \displaystyle\frac{k_i(k_i-1)}{n(n+1)},&i=j.
    \end{cases}
\end{equation}
Lemma~\ref{lem:q-probability} identifies this kernel as the winner--loser type distribution for two individuals sampled uniformly without replacement from a population with degree vector $\kvec$.
The \emph{selection} operator $S_n:\V_n\to\V_{n+1}$ is
\begin{equation}\label{eq:S-operator}
    (S_nv)(\kvec):=\sum_{\substack{j\in[d]:k_j\geq1}}\sum_{i=1}^d
    q_{ij}^{(n+1)}(\kvec)v(\kvec-\e_j),
    \qquad \kvec\in\Kc_{n+1}.
\end{equation}
The restriction $k_j\geq1$ ensures that every argument of $v$ belongs to $\Kc_n$.
Note that $q_{ij}^{(n+1)}(\kvec)=0$ whenever $k_j=0$.
For a fixed colouring $\ivec=(i_1,\ldots,i_{n+1})$ with $\operatorname{deg}(\ivec)=\kvec$, let $(U_1,U_2)$ be an ordered uniform sample of distinct labels from $[n+1]$.
Thus $i_{U_h}$ is the type of the sampled individual $U_h$, $h=1,2$.
Lemma~\ref{lem:q-probability} gives
\[
    (S_nv)(\kvec)=\E[v(\kvec-\e_{\ell(i_{U_1},i_{U_2})})].
\]

For $n\geq2$, define the \emph{coalescent} operator $C_n:\V_n\to\V_{n-1}$ by
\begin{equation}\label{eq:C-operator}
    (C_nv)(\kvec):=\sum_{j=1}^d\frac{k_j}{n-1}v(\kvec+\e_j),
    \qquad \kvec\in\Kc_{n-1}.
\end{equation}
For a fixed colouring $\ivec\in[d]^{n-1}$ with degree vector $\kvec$ and a uniform label $U_1\in[n-1]$,
\[
    (C_nv)(\kvec)=\E[v(\kvec+\e_{i_{U_1}})],
\]
since $\Pp(i_{U_1}=j)=k_j/(n-1)$.
Both operators are linear and satisfy
\begin{equation}\label{eqnDomainsOfOperators}
    d(S_nv)=n+1,
    \qquad d(C_nv)=n-1\quad(n\geq2).
\end{equation}
Here the index $d(v)$ records the Bernstein representation, not necessarily the minimal polynomial degree of $H(v,\cdot)$ after restriction to the simplex.

\begin{definition}[Multivariate Bernstein coefficient process]\label{def:coefficient-process}
Fix an integer $d\geq2$, a tournament matrix $T$, parameters $\lambda,\gamma\in[0,\infty)$, an integer $m\in\N$, and $v_0\in\V_m$.
The \emph{multivariate Bernstein coefficient process} is the pure-jump process $V=(V_t)_{t\geq0}$ on $\V$, started from $V_0=v_0$, whose extended jump operator on bounded measurable $F:\V\to\R$ is
\begin{equation}\label{eq:V-generator}
    \mathcal BF(v):=\lambda n[F(S_nv)-F(v)]
       +\gamma\binom n2[F(C_nv)-F(v)],
    \qquad v\in\V_n,\quad n\in\N,
\end{equation}
with the coalescence term omitted when $n=1$.
\end{definition}

Write $L_t:=d(V_t)$, so $L_0=m$.
The \emph{degree process} $L$ has transitions
\begin{equation}\label{eq:degree-rates}
n\longmapsto n+1\quad\text{at rate }\lambda n,
\qquad
n\longmapsto n-1\quad\text{at rate }\gamma\binom n2,
\end{equation}
where the second transition is present only for $n\geq2$.
From degree $n$, the total jump rate is $q_n:=\lambda n+\gamma\binom n2$.
If $q_n>0$, wait an exponential time of rate $q_n$ and apply $S_n$ or $C_n$ with probabilities $\lambda n/q_n$ and $\gamma\binom n2/q_n$, respectively, using independent holding-time and choice variables at successive jumps.
If $q_n=0$, the state is absorbing.
Lemma~\ref{lem:degree-nonexplosion} proves that this construction is nonexplosive.

\begin{figure}[tbp]
\centering
\reviewfigure{.6\textwidth}{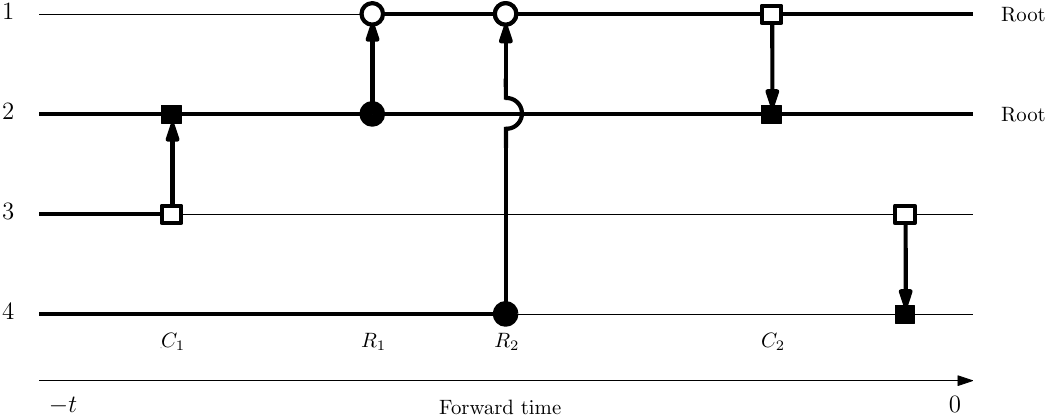}
\caption{Graphical representation of an uncoloured finite-population  ancestral
		selection graph for two sampled individuals (roots 1 and 2).
		Squares indicate competitions, with the white square marking the challenger
		$U_1$ and the black square marking the target $U_2$ (whose type is determined by the tournament outcome).
		Circles indicate neutral reproductions, with the black circle marking the
		reproducing individual $U_1$ and the white circle marking the target $U_2$  (whose type will be copied from $U_1$).
		Bold lines show the ancestral lines that are followed for the sample.
		At $C_1$, individuals $U_1=3$ and $U_2=2$ compete.
		At $R_1$, individual $U_1=2$ reproduces and replaces individual $U_2=1$.
		At $R_2$, the event involves individuals 4 and 1, showing that a reproduction event
		may involve an individual outside the sampled roots.
		The additional pair of squares near the right represents a competition between individuals 3 and 4 outside the ancestry of the sample; hence, the corresponding right-hand segments are not drawn in bold.
		Note that assigning types to the ancestral lines shown in bold at the left uniquely determines the types of the sampled individuals.}
\label{fig:asg-graph}
\end{figure}

\begin{figure}[tbp]
\centering
\reviewfigure{.6\textwidth}{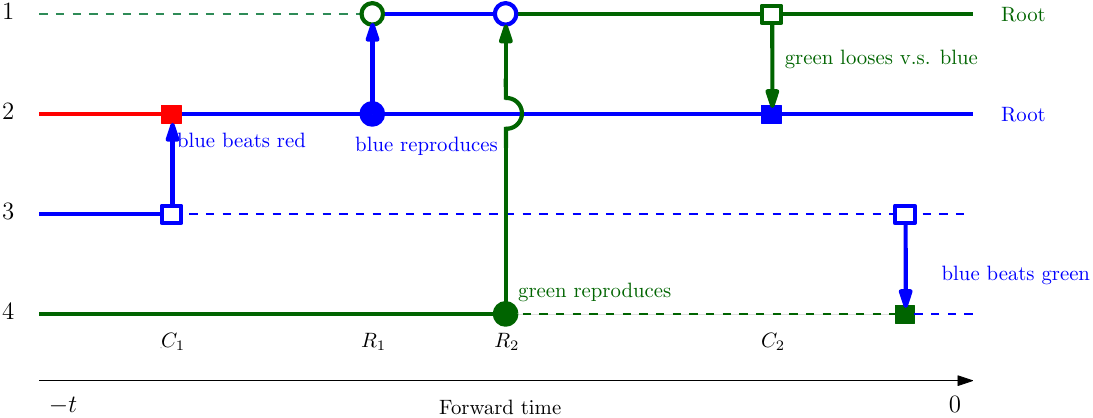}
\caption{Evolution of types along the same realization as in
		Figure~\ref{fig:asg-graph}.
		Tournament used here:
		blue beats green and red, red beats green.
		At the left, individuals 1 and 4 are green, individual 2 is red, individual 3 is blue.
		The arrows are coloured by the type that challenges or reproduces (the $U_1$'s).
		At $C_1$, individual 3 has type blue and individual 2 has type red, so
		\emph{blue beats red}; since the target is $U_2=2$, individual 2 changes to blue.
		At $R_1$, individual 2 is blue and reproduces to individual 1, so individual 1
		changes from green to blue.
		At $R_2$, individual 4 is green and reproduces to individual 1, so individual 1
		changes back to green.
		At $C_2$, individual 1 is green and individual 2 is blue, since 
		\emph{green loses to blue}; the target $U_2=2$ remains blue.
		Finally, in the competition outside the sample, individual 3 challenges 
		individual 4 and wins; this changes individual 4 to blue,
		but it does not affect the sampled roots.}
\label{fig:asg-types}
\end{figure}

\begin{remark}[Ancestral interpretation]
Backward in time, the degree change under selection represents the addition of one potential parental line.
When types are propagated towards the sample, the two parental types are replaced by their tournament winner, so $S_n$ averages over deletion of the losing type.
A coalescence instead merges two lines backward in time; read forward, the type of their common parent is duplicated, as described by $C_n$.
Figures~\ref{fig:asg-graph} and \ref{fig:asg-types} illustrate the corresponding mechanisms in a finite population, not a realization of the limiting degree process.
In the finite graph, the challenger may already be among the traced lines, so a competition need not increase their number.
\end{remark}

\subsection{Main result}\label{subsecCoefficientDynamics}

Recall that a continuous process $\X$ with values in $\Delta_d$
solves the martingale problem for
$(\mathcal L,C^2(\Delta_d))$ with initial state
$\x\in\Delta_d$ if $\X_0=\x$ almost surely and,
for every $f\in C^2(\Delta_d)$, the process
\begin{equation}\label{eq:forward-martingale-definition}
f(\X_t)-f(\x)-\int_0^t\mathcal Lf(\X_s)\,\dd s,
\qquad t\geq0,
\end{equation}
is a martingale with respect to the natural filtration
of $\X$.
Existence of such a process is established in
Appendix~\ref{app:forward}.

Our main result proves the multivariate Bernstein duality between the frequency process and the multivariate Bernstein coefficient process.

\begin{theorem}[Multivariate Bernstein duality]\label{thm:bernstein-duality}
Let $\X$ be any continuous solution of the martingale problem for $(\mathcal L,C^2(\Delta_d))$, and let $V$ be the coefficient process of Definition~\ref{def:coefficient-process}, with the same tournament and parameters $\lambda,\gamma\in[0,\infty)$.
Then, for every $t\geq0$, $n\in\N$, $v\in\V_n$, and $\x\in\Delta_d$,
\begin{equation}\label{eq:semigroup-duality}
    \E_{\x}[H(v,\X_t)]=\E_v[H(V_t,\x)].
\end{equation}
More explicitly,
\[
    \E_{\x}\Bigl[\sum_{\kvec\in\Kc_n}
    v(\kvec)\binomulti n{\kvec}\X_t^{\kvec}\Bigr]
=\E_v\Bigl[\sum_{\kvec\in\Kc_{L_t}}
V_t(\kvec)\binomulti {L_t}{\kvec}\x^{\kvec}\Bigr].
\]
\end{theorem}

The two expectations refer to different processes with deterministic initial states $\X_0=\x$ and $V_0=v$, respectively.

\section{Multinomial sampling and operator identities}
\label{sec:bernstein}

In this section, we establish probabilistic sampling identities
that relate independent multinomial sampling to sampling
without replacement.
These identities provide probabilistic interpretations of
the selection and coalescent operators.
Although the arguments are elementary, we include complete proofs for clarity.
The main results are
Propositions~\ref{prop:selection-intertwining}
and~\ref{prop:coalescence-intertwining},
from which one trivially establishes the generator identities underlying
the multivariate Bernstein duality.

\subsection{Elementary multinomial identities}

The next lemma records the derivatives used in the generator calculation.
Differentiating corresponds to fixing one or two sample types and averaging over the remaining independent types.
The derivatives are taken in the ambient space $\R^d$ using the polynomial extension in \eqref{eq:H-definition}.

\begin{lemma}[Multinomial derivative identities]\label{lem:derivative-identities}
Let $n\geq1$ and $v\in\V_n$.
For every $i,j\in[d]$ and $\x\in\Delta_d$,
\begin{align}
    \partial_jH(v,\x)
    &=n\E[v(\Mvec_{n-1}(\x)+\e_j)],
    \label{eq:first-derivative-H}\\
    \partial_{ij}H(v,\x)
    &=n(n-1)\E[v(\Mvec_{n-2}(\x)+\e_i+\e_j)],
    \label{eq:second-derivative-H}
\end{align}
where \eqref{eq:second-derivative-H} is used only when $n\geq2$.
\end{lemma}

\begin{proof}
Differentiate the finite sum in \eqref{eq:H-definition} term by term.
For the first derivative, only indices with positive $j$th coordinate contribute.
Write each such index as $\kvec+\e_j$, with $\kvec\in\Kc_{n-1}$.
Then
\[
    (k_j+1)\binomulti n{\kvec+\e_j}
    =(k_j+1)\frac{n!}{(k_j+1)!\prod_{i\neq j}k_i!}
    =n\binomulti{n-1}{\kvec}.
\]
Consequently,
\[
    \partial_jH(v,\x)
    =n\sum_{\kvec\in\Kc_{n-1}}
    v(\kvec+\e_j)\binomulti{n-1}{\kvec}\x^{\kvec},
\]
which is \eqref{eq:first-derivative-H}.
Applying the same calculation once more gives \eqref{eq:second-derivative-H}.
This also covers $i=j$, because differentiating twice multiplies a monomial by the corresponding falling factorial.
\end{proof}

For $1\leq\ell\leq m$, write $\Uvec^{(\ell)}=(U_1,\ldots,U_\ell)$ for an ordered uniform sample from $[m]$ without replacement, that is,
\begin{equation}\label{eq:ordered-label-law}
    \Pp(\Uvec^{(\ell)}=\uvec)=\frac{1}{(m)_{\underline\ell}},
    \qquad \uvec=(u_1,\ldots,u_\ell)\in[m]^\ell,
    \quad u_r\neq u_s\ \text{for }r\neq s,
\end{equation}
where $(m)_{\underline\ell}:=m(m-1)\cdots(m-\ell+1)$.
For simplicity on exposition, we suppress the dependence on $m$ in this notation.

The following lemma is the key ingredient in the proof of the multivariate Bernstein duality, Theorem~\ref{thm:bernstein-duality}.
It shows that, if $m$ individuals are coloured independently with common law $\x$, then the types of $\ell$ individuals sampled uniformly without replacement are still independent with law $\x$.
Moreover, the degree vector of the remaining $m-\ell$ individuals is independent of the sampled types and has a multinomial distribution.

\begin{lemma}[Exchangeable marking]\label{lem:exchangeable-marking}
Let $1\leq\ell\leq m$, let $\x\in\Delta_d$, and let $I_1,\ldots,I_m$ be independent random types with law $\x$.
Independently of $(I_1,\ldots,I_m)$, let $\Uvec^{(\ell)}=(U_1,\ldots,U_\ell)$ be as above.
The marked individuals $U_1,\ldots,U_\ell$ have types $I_{U_1},\ldots,I_{U_\ell}$.
Define the degree vector of the unmarked individuals by
$\Mvec:=\sum_{u\in[m]\setminus\{U_1,\ldots,U_\ell\}}\e_{I_u}$.
Then the marked types are independent with common law $\x$, and $\Mvec$ has law $\Mult(m-\ell,\x)$ independently of all the marked types.
Equivalently,
\begin{equation}\label{eq:marking-distributional}
    (I_{U_1},\ldots,I_{U_\ell},\Mvec)
    \eqdist
    (I_1,\ldots,I_\ell,\Mvec_{m-\ell}(\x)),
\end{equation}
where $\Mvec_{m-\ell}(\x)$ is independent of $I_1,\ldots,I_m$.
\end{lemma}

\begin{proof}
Fix $\uvec=(u_1,\ldots,u_\ell)\in[m]^\ell$ with all entries distinct, $(i_1,\ldots,i_\ell)\in[d]^\ell$, and $\kvec\in\Kc_{m-\ell}$.
Conditional on $\Uvec^{(\ell)}=\uvec$, the marked types are $I_{u_1},\ldots,I_{u_\ell}$ and
$\Mvec=\sum_{u\in[m]\setminus\{u_1,\ldots,u_\ell\}}\e_{I_u}$.
This conditioning does not alter the law of the type vector because the sample of labels is independent of it.
The marked and unmarked types now belong to disjoint, deterministic subfamilies of the independent family $(I_u)_{u\in[m]}$.
Hence, by \eqref{eq:multinomial-realization},
\begin{equation}\label{eq:marking-joint-factorization}
\begin{aligned}
&\Pp\bigl(\Uvec^{(\ell)}=\uvec,
I_{U_1}=i_1,\ldots,I_{U_\ell}=i_\ell,\Mvec=\kvec\bigr)=\frac1{(m)_{\underline\ell}}
\left(\prod_{h=1}^\ell x_{i_h}\right)
\binomulti{m-\ell}{\kvec}\x^{\kvec}.
\end{aligned}
\end{equation}
Summing over the $(m)_{\underline\ell}$ possible label samples gives
\[
    \Pp(I_{U_1}=i_1,\ldots,I_{U_\ell}=i_\ell,\Mvec=\kvec)
    =\left(\prod_{h=1}^\ell x_{i_h}\right)
      \Pp(\Mvec_{m-\ell}(\x)=\kvec).
\]
This factorization proves the asserted independence and all the marginal laws, and therefore \eqref{eq:marking-distributional}.
It also covers $\ell=m$, when the unmarked degree vector is $\Mvec=\zerovec$.
\end{proof}

\subsection{The winner--loser selection operator}

Recall the winner--loser maps and kernel given in \eqref{eq:winner-loser-map}--\eqref{eq:q-winner-loser}.
For a fixed population composition, the next lemma counts the possible winning and losing pairs; their arrangement among labels does not matter.

\begin{lemma}[Conditional winner--loser law]\label{lem:q-probability}
Fix $n\geq1$, $\kvec\in\Kc_{n+1}$, and a deterministic type assignment $\ivec=(i_1,\ldots,i_{n+1})\in[d]^{n+1}$ with $\operatorname{deg}(\ivec)=\kvec$.
Choose $\Uvec^{(2)}=(U_1,U_2)$ as in \eqref{eq:ordered-label-law} with $m=n+1$, and set
$\mathsf W:=w(i_{U_1},i_{U_2})$ and $\mathsf L:=\ell(i_{U_1},i_{U_2})$.
Thus $\mathsf W$ is the retained type and $\mathsf L$ is the removed type.
Then, for all $i,j\in[d]$,
\begin{equation}\label{eq:fixed-population-winner-law}
    \Pp(\mathsf W=i,\mathsf L=j)=q_{ij}^{(n+1)}(\kvec).
\end{equation}
In particular, $q^{(n+1)}(\kvec)$ is a probability distribution that depends only on the degree vector, not on the arrangement of the types among individual labels.
For a random type assignment and an independent uniform sample $\Uvec^{(2)}$, the same kernel gives the conditional law of $(\mathsf W,\mathsf L)$ given the degree vector.
\end{lemma}

\begin{proof}
There are $n(n+1)$ ordered pairs of distinct individual labels, all equally likely.
For $i\neq j$, there are $k_ik_j$ pairs with ordered types $(i,j)$ and another $k_ik_j$ pairs with ordered types $(j,i)$.
For $i=j$, there are $k_i(k_i-1)$ pairs with ordered types $(i,i)$.
Consequently,
\begin{equation}\label{eq:conditional-pair-types}
    \Pp(i_{U_1}=i,i_{U_2}=j)
    =\frac{k_i(k_j-\delta_{ij})}{n(n+1)}.
\end{equation}
For distinct $i,j$, the event $\{\mathsf W=i,\mathsf L=j\}$ is empty if $T_{ij}=0$.
If $T_{ij}=1$, then
\[
    \{\mathsf W=i,\mathsf L=j\}
    =\{i_{U_1}=i,i_{U_2}=j\}
      \cup\{i_{U_1}=j,i_{U_2}=i\}.
\]
These two events are disjoint because $i\neq j$, so their probabilities add.
Equation~\eqref{eq:conditional-pair-types} therefore gives $2k_ik_jT_{ij}/[n(n+1)]$.
For $i=j$, the event $\{\mathsf W=i,\mathsf L=i\}$ is precisely $\{i_{U_1}=i,i_{U_2}=i\}$ and has probability $k_i(k_i-1)/[n(n+1)]$.
This proves \eqref{eq:fixed-population-winner-law}.
The entries are nonnegative, and
\begin{align*}
    n(n+1)\sum_{i,j=1}^d q_{ij}^{(n+1)}(\kvec)
    &=\sum_{j=1}^d k_j(k_j-1)
      +2\sum_{i<j}k_i k_j(T_{ij}+T_{ji})\\
    &=\left(\sum_{j=1}^d k_j\right)^2-\sum_{j=1}^d k_j
      =n(n+1).
\end{align*}
For a random assignment, condition first on the full assignment.
Independence of the sampled labels makes the preceding computation applicable, and its result is a function of the degree vector alone.
Conditioning further on that vector proves the last assertion by the tower property.
Conditional probabilities at degree vectors of probability zero may be assigned the same kernel; all conditional identities are asserted almost surely.
\end{proof}

The next lemma applies the preceding fixed-colouring calculation to independent random colours with law $\x$.
After removing the loser of \emph{two distinct individuals sampled uniformly}, the resulting degree vector has the same law as an independent sample of size $n-1$ together with the winner of two additional independent \emph{types}.

\begin{lemma}[Winner--loser sampling identity]\label{lem:winner-loser-identity}
Fix $n\geq1$ and $\x\in\Delta_d$.
Let $I_1,\ldots,I_{n+1}$ be independent random types with law $\x$ and set $\Mvec_{n+1}(\x):=\sum_{u=1}^{n+1}\e_{I_u}$.
Independently of this type vector, choose an ordered uniform sample $\Uvec^{(2)}=(U_1,U_2)$ as in \eqref{eq:ordered-label-law} with $m=n+1$, and define $\mathsf W:=w(I_{U_1},I_{U_2})$ and $\mathsf L:=\ell(I_{U_1},I_{U_2})$.
Let $\Mvec_{n-1}(\x)$ be a further multinomial vector independent of $I_1,\ldots,I_{n+1}$.
Then
\begin{equation}\label{eq:winner-loser-distributional}
    (\mathsf W,\mathsf L,\Mvec_{n+1}(\x)-\e_{\mathsf L})
    \eqdist
    \bigl(w(I_1,I_2),\ell(I_1,I_2),
          \Mvec_{n-1}(\x)+\e_{w(I_1,I_2)}\bigr).
\end{equation}
In particular, for every $v\in\V_n$,
\begin{equation}\label{eq:HS-winner-representation}
    H(S_nv,\x)
    =\E[v(\Mvec_{n-1}(\x)+\e_{w(I_1,I_2)})].
\end{equation}
\end{lemma}

\begin{proof}
Define the degree vector of the unmarked individuals by $\Mvec:=\sum_{u\in[n+1]\setminus\{U_1,U_2\}}\e_{I_u}$.
Lemma~\ref{lem:exchangeable-marking}, with $m=n+1$ and $\ell=2$, gives
\begin{equation}\label{eqnTypesOf2SampleAreiid}
    (I_{U_1},I_{U_2},\Mvec)
    \eqdist(I_1,I_2,\Mvec_{n-1}(\x)),
\end{equation}
where the three components on the right are independent.
By construction, $\Mvec_{n+1}(\x)=\Mvec+\e_{I_{U_1}}+\e_{I_{U_2}}$.
For every pair of types $i,j$, we have $\e_i+\e_j=\e_{w(i,j)}+\e_{\ell(i,j)}$.
For distinct types this just reorders the two summands, and for equal types both sides equal $2\e_i$.
It follows that
\begin{equation}\label{eq:deleting-loser-pathwise}
    \Mvec_{n+1}(\x)-\e_{\mathsf L}
    =\Mvec+\e_{w(I_{U_1},I_{U_2})}
    =\Mvec+\e_{\mathsf W}
    \qquad\text{almost surely}.
\end{equation}
Applying the deterministic map $(i,j,\kvec)\mapsto(w(i,j),\ell(i,j),\kvec+\e_{w(i,j)})$ to \eqref{eqnTypesOf2SampleAreiid}, and using \eqref{eq:deleting-loser-pathwise}, proves \eqref{eq:winner-loser-distributional}.

To obtain the polynomial identity, put $\mathcal F_I:=\sigma(I_1,\ldots,I_{n+1})$.
Conditional on this sigma-field, every ordered pair of distinct individual labels has probability $1/[n(n+1)]$, because $\Uvec^{(2)}$ is independent of $\mathcal F_I$.
Lemma~\ref{lem:q-probability} therefore yields
$\Pp(\mathsf W=i,\mathsf L=j\mid\mathcal F_I)=q_{ij}^{(n+1)}(\Mvec_{n+1}(\x))$.
The right-hand side is measurable with respect to $\sigma(\Mvec_{n+1}(\x))$.
Thus, the tower property gives
$\Pp(\mathsf W=i,\mathsf L=j\mid\Mvec_{n+1}(\x))=q_{ij}^{(n+1)}(\Mvec_{n+1}(\x))$ almost surely.
For every $\kvec$ with $\Pp(\Mvec_{n+1}(\x)=\kvec)>0$, it follows that
\begin{equation}\label{eq:selection-conditional-expectation}
\begin{aligned}
&\E[v(\Mvec_{n+1}(\x)-\e_{\mathsf L})\mid\Mvec_{n+1}(\x)=\kvec]\\
&\qquad=\sum_{\substack{j\in[d]:k_j\geq1}}\sum_{i=1}^d
q_{ij}^{(n+1)}(\kvec)v(\kvec-\e_j)=(S_nv)(\kvec).
\end{aligned}
\end{equation}
All variables are integrable because $v$ has finitely many values.
Taking expectations and using \eqref{eq:H-definition} gives
\begin{align*}
H(S_nv,\x)
&=\E[(S_nv)(\Mvec_{n+1}(\x))]\\
&=\E[v(\Mvec_{n+1}(\x)-\e_{\mathsf L})]\\
&=\E[v(\Mvec_{n-1}(\x)+\e_{w(I_1,I_2)})],
\end{align*}
where the last equality follows from \eqref{eq:winner-loser-distributional}.
This proves \eqref{eq:HS-winner-representation}.
For $n=1$, both unmarked and reference degree vectors are identically $\zerovec$, and the same proof applies.
\end{proof}

The following calculation shows that ancestral selection produces exactly the selection part of the forward generator.

\begin{proposition}[Selection intertwining identity]\label{prop:selection-intertwining}
For every $n\geq1$, $v\in\V_n$, and $\x\in\Delta_d$,
\begin{equation}\label{eq:selection-polynomial-identity}
    H(S_nv,\x)-H(v,\x)
    =\sum_{i=1}^d\sum_{j\neq i}
    x_ix_jA_{ij}\E[v(\Mvec_{n-1}(\x)+\e_i)].
\end{equation}
Consequently,
\begin{equation}\label{eq:selection-generator-match}
    \lambda\sum_{j=1}^dx_j(A\x)_j\partial_jH(v,\x)
    =\lambda n[H(S_nv,\x)-H(v,\x)].
\end{equation}
\end{proposition}

\begin{proof}
For $j\in[d]$, set $g_j:=\E[v(\Mvec_{n-1}(\x)+\e_j)]$.
By \eqref{eq:HS-winner-representation}, with the reference multinomial vector independent of $I_1,I_2$,
\begin{align*}
H(S_nv,\x)
&=\sum_{j=1}^d x_j^2g_j
 +\sum_{\substack{i,j\in[d]\\i\neq j}}x_ix_jg_{w(i,j)}=\sum_{j=1}^d\left(x_j^2+2\sum_{i\neq j}x_jx_iT_{ji}\right)g_j.
\end{align*}
The second equality groups the ordered type pairs according to their common winner $j$.
On the other hand, independence implies $\Mvec_{n-1}(\x)+\e_{I_1}\eqdist\Mvec_n(\x)$, so $H(v,\x)=\sum_jx_jg_j$.
Since
\[
x_j=x_j^2+\sum_{i\neq j}x_jx_i,
\qquad
2T_{ji}-1=T_{ji}-T_{ij}=A_{ji}\quad(i\neq j),
\]
subtraction gives \eqref{eq:selection-polynomial-identity}.
Finally, $\partial_jH(v,\x)=ng_j$ by Lemma~\ref{lem:derivative-identities} and $(A\x)_j=\sum_iA_{ji}x_i$, which gives \eqref{eq:selection-generator-match}.
\end{proof}

\subsection{The Wright--Fisher coalescent operator}

We use the operator $C_n$ from \eqref{eq:C-operator}.
Recall that for $n\geq2$, a binary coalescence read forward corresponds to selecting one of the $n-1$ parental leaves uniformly and duplicating its type.
The operator averages over the type of that leaf.
More precisely, fix $\ivec\in[d]^{n-1}$ with $\operatorname{deg}(\ivec)=\kvec$ and let $\Uvec^{(1)}=U_1$ be a uniform one-label sample from $[n-1]$.
Then $\Pp(i_{U_1}=j)=k_j/(n-1)$ and $\e_{i_{U_1}}\eqdist\Mvec_1(\kvec/(n-1))$.
Consequently,
\begin{equation}\label{eq:coalescence-fixed-expectation}
    (C_nv)(\kvec)=\E[v(\kvec+\e_{i_{U_1}})]
    =\E\left[v\left(\kvec+\Mvec_1\left(\frac{\kvec}{n-1}\right)\right)\right].
\end{equation}
The last expression is exactly \eqref{eq:C-operator}.
Here $\kvec/(n-1)\in\Delta_d$ is the empirical frequency vector of the fixed population.
When that population is itself sampled from $\x$, the identity is conditional on its degree vector, as shown next.
The next lemma states that the degree vector \emph{after duplication} has the same law as an independent sample of size $n-2$ together with two copies of an independently sampled type.

\begin{lemma}[Duplication sampling identity]\label{lem:duplication-identity}
Let $n\geq2$ and $\x\in\Delta_d$.
Take independent random types $I_1,\ldots,I_{n-1}$ with law $\x$ and set $\Mvec_{n-1}(\x):=\sum_{u=1}^{n-1}\e_{I_u}$.
Independently of the full type vector, choose a uniform one-label sample $\Uvec^{(1)}=U_1$ from $[n-1]$.
Then
\begin{equation}\label{eq:duplication-distributional}
    \Mvec_{n-1}(\x)+\e_{I_{U_1}}
    \eqdist\Mvec_{n-2}(\x)+2\e_{I_1},
\end{equation}
where the reference multinomial vector on the right is chosen independently of $I_1,\ldots,I_{n-1}$.
Consequently, for $v\in\V_n$,
\begin{equation}\label{eq:HC-duplication}
    H(C_nv,\x)=\sum_{j=1}^dx_j\E[v(\Mvec_{n-2}(\x)+2\e_j)].
\end{equation}
\end{lemma}

\begin{proof}
Put $\Mvec:=\sum_{u\in[n-1]\setminus\{U_1\}}\e_{I_u}$.
Lemma~\ref{lem:exchangeable-marking}, with $m=n-1$ and $\ell=1$, gives
\[
    (I_{U_1},\Mvec)\eqdist(I_1,\Mvec_{n-2}(\x)),
\]
with independent components on the right.
The pathwise equality $\Mvec_{n-1}(\x)+\e_{I_{U_1}}=\Mvec+2\e_{I_{U_1}}$ therefore proves \eqref{eq:duplication-distributional}.
Conditional on the full type assignment, a uniform individual label has type $j$ with probability $[\Mvec_{n-1}(\x)]_j/(n-1)$.
The tower property implies, for every $\kvec$ of positive multinomial probability,
$\Pp(I_{U_1}=j\mid\Mvec_{n-1}(\x)=\kvec)=\frac{k_j}{n-1}$.
Hence, by \eqref{eq:C-operator},
\begin{equation}\label{eq:coalescence-conditional-expectation}
\begin{aligned}
&\E[v(\Mvec_{n-1}(\x)+\e_{I_{U_1}})\mid\Mvec_{n-1}(\x)=\kvec]=\sum_{j=1}^d\frac{k_j}{n-1}v(\kvec+\e_j)=(C_nv)(\kvec).
\end{aligned}
\end{equation}
Since $C_nv\in\V_{n-1}$, equation~\eqref{eq:H-definition} gives
\begin{equation}\label{eq:H-coalescence-domain}
    H(C_nv,\x)
    =\sum_{\kvec\in\Kc_{n-1}}(C_nv)(\kvec)\binomulti{n-1}{\kvec}\x^{\kvec}
    =\E[(C_nv)(\Mvec_{n-1}(\x))].
\end{equation}
Taking expectations in \eqref{eq:coalescence-conditional-expectation}, applying \eqref{eq:duplication-distributional}, and conditioning on the independent reference type $I_1$, we obtain
\begin{align*}
    H(C_nv,\x)
    &=\E[(C_nv)(\Mvec_{n-1}(\x))]\\
    &=\E[v(\Mvec_{n-1}(\x)+\e_{I_{U_1}})]\\
    &=\E[v(\Mvec_{n-2}(\x)+2\e_{I_1})]\\
    &=\sum_{j=1}^d x_j\E[v(\Mvec_{n-2}(\x)+2\e_j)].
\end{align*}
For $n=2$, the unmarked and reference multinomial vectors are zero, so the calculation remains valid.
\end{proof}

The next result shows that duplicating a common parent reproduces the Wright--Fisher covariance term.

\begin{proposition}[Coalescent intertwining identity]\label{prop:coalescence-intertwining}
For every $n\geq2$, $v\in\V_n$, and $\x\in\Delta_d$,
\begin{equation}\label{eq:coalescence-generator-match}
    \frac{\gamma}{2}\sum_{i,j=1}^d
    x_i(\delta_{ij}-x_j)\partial_{ij}H(v,\x)
    =\gamma\binom{n}{2}[H(C_nv,\x)-H(v,\x)].
\end{equation}
\end{proposition}

\begin{proof}
By \eqref{eq:second-derivative-H}, the left-hand side is
\[
\gamma\binom{n}{2}\Bigg[
\sum_{j=1}^dx_j\E[v(\Mvec_{n-2}(\x)+2\e_j)]
-\sum_{i,j=1}^dx_ix_j\E[v(\Mvec_{n-2}(\x)+\e_i+\e_j)]
\Bigg].
\]
The first term in brackets is $H(C_nv,\x)$ by Lemma~\ref{lem:duplication-identity}.
For the second, take independent types $I_1,I_2$ with law $\x$, independent of $\Mvec_{n-2}(\x)$.
Then $\Mvec_{n-2}(\x)+\e_{I_1}+\e_{I_2}\eqdist\Mvec_n(\x)$, so the second term equals $H(v,\x)$.
This proves \eqref{eq:coalescence-generator-match}.
\end{proof}

\subsection{Stability of the coefficient operators}

We now show that both coefficient operators preserve the convex hull of the input values. 

\begin{lemma}[Contraction and vertex values]\label{lem:operator-contraction}
For every $v\in\V_n$,
\begin{equation}\label{eq:operator-contraction}
    \norm{S_nv}_\infty\leq\norm{v}_\infty,
    \qquad
    \norm{C_nv}_\infty\leq\norm{v}_\infty,
\end{equation}
where the second inequality is used for $n\geq2$.
Both operators take values in the convex hull of the values of $v$.
Moreover, for every $j\in[d]$,
\begin{equation}\label{eq:operator-vertex-values}
    (S_nv)((n+1)\e_j)=v(n\e_j),
    \qquad
    (C_nv)((n-1)\e_j)=v(n\e_j),
\end{equation}
again using the second identity only for $n\geq2$.
\end{lemma}

\begin{proof}
Put $m_v:=\min_{\hvec\in\Kc_n}v(\hvec)$ and $M_v:=\max_{\hvec\in\Kc_n}v(\hvec)$.
Fix $\kvec\in\Kc_{n+1}$.
Every index in \eqref{eq:S-operator} has $k_j\geq1$, so $\kvec-\e_j$ has nonnegative integer coordinates summing to $n$.
Thus $\kvec-\e_j\in\Kc_n$ and $|v(\kvec-\e_j)|\leq\norm v_\infty$.
Using nonnegative weights and their sum, which is one by Lemma~\ref{lem:q-probability}, gives
\begin{align*}
|(S_nv)(\kvec)|
&\leq\sum_{j:k_j\geq1}\sum_{i=1}^d
q_{ij}^{(n+1)}(\kvec)|v(\kvec-\e_j)|\\
&\leq\norm v_\infty\sum_{j:k_j\geq1}\sum_{i=1}^d
q_{ij}^{(n+1)}(\kvec)=\norm v_\infty.
\end{align*}
Since the bound is independent of $\kvec$, this proves the first inequality in \eqref{eq:operator-contraction}.
Similarly, for $n\geq2$ and $\kvec\in\Kc_{n-1}$, every vector $\kvec+\e_j$ belongs to $\Kc_n$ and $\sum_jk_j/(n-1)=1$.
Therefore
\[
    |(C_nv)(\kvec)|\leq\sum_{j=1}^d\frac{k_j}{n-1}|v(\kvec+\e_j)|
    \leq\norm v_\infty.
\]
Replacing the absolute-value bound by $m_v\leq v(\hvec)\leq M_v$ in the same computations proves the convex-hull assertion.
At $(n+1)\e_j$, the winner and loser types are necessarily $(j,j)$, so $q_{jj}^{(n+1)}((n+1)\e_j)=1$ (recall \eqref{eq:q-winner-loser}).
At $(n-1)\e_j$, the duplicated leaf necessarily has type $j$, so $k_j/(n-1)=1$ in \eqref{eq:C-operator}, and all other weights are zero.
These observations give \eqref{eq:operator-vertex-values}.
\end{proof}

\subsection{Nonexplosion and moments of the coefficient process}\label{sec:coefficient}

The proof of Theorem~\ref{thm:bernstein-duality} requires nonexplosion and an integrable bound for $\sup_{s\leq T}L_s^2$.
We use the same topological convention as \cite[Remark~4.6]{CorderoHummelSchertzer2022}: retain the degree and use the usual finite-dimensional topology within each coefficient space.
For each $n\in\N$, fix an enumeration $\kvec^{(n,1)},\ldots,\kvec^{(n,D_n)}$ of $\Kc_n$, where $D_n=\binom{n+d-1}{d-1}$.
For $v\in\V_n$, let $\widetilde v\in\R^{\N}$ have entries $\widetilde v_r=v(\kvec^{(n,r)})$ for $r\leq D_n$ and $\widetilde v_r=0$ otherwise.
Equip $\V$ with the metric
\begin{equation}\label{eq:V-metric}
    d_{\V}(v,w):=\|\widetilde v-\widetilde w\|_\infty+|d(v)-d(w)|.
\end{equation}
The degree term prevents different coefficient domains from being identified by zero padding.
Note that $d_{\V}$ is the usual sup-norm metric on every $\V_n$, and different degrees are at least distance one apart.
Thus it induces the same disjoint-union topology as the coefficient-array construction in the cited remark.
A Cauchy sequence has eventually constant degree and then converges in a finite-dimensional complete space.
Coefficient arrays with rational entries form a countable dense subset, so $\V$ is complete and separable.
In particular, a function on $\V$ is continuous exactly when its restriction to each $\V_n$ is continuous.

Recall the processes $V$ and $L$ from Definition~\ref{def:coefficient-process}.
The following result shows that a finite ancestral sample has only finitely many events on each bounded time interval, and its maximum size has finite moments.

\begin{lemma}[Nonexplosion and degree moments]\label{lem:degree-nonexplosion}
The processes $L$ and $V$ are nonexplosive.
Moreover, for every initial degree $n\in\N$, $T>0$, integer $p\geq1$, and $\lambda,\gamma\in[0,\infty)$,
\begin{equation}\label{eq:degree-moment-bound}
    \E_n\left[\sup_{0\leq s\leq T}L_s^p\right]<\infty.
\end{equation}
\end{lemma}

\begin{proof}
If $\lambda=0$, then $L_t\leq n$ and at most $n-1$ coalescences can occur, so the claims follow immediately.
Assume $\lambda>0$ and let $Y$ be the nonexplosive Yule process with $Y_0=n$ and birth rate $\lambda y$ in state $y$.
Then $\E_n[Y_T^p]<\infty$ for every finite $T$ and positive integer $p$, since $Y_T-n$ has a negative binomial distribution (see \cite[Chapter~6, Section~6.4, Example~6.8]{Ross2010}).
One can couple $(L,Y)$, initially before any possible explosion, so that at each positive jump of $L$, there is a positive jump of $Y$.
It follows that $L_s\leq Y_s\leq Y_T$ for $0\leq s\leq T$. 
This implies that $L$ has only finitely many births on $[0,T]$, and $\E_n\left[\sup_{0\leq s\leq T}L_s^p\right]
\leq\E_n[Y_T^p]<\infty$.

\end{proof}

For a bounded measurable $F$, the right-hand side of \eqref{eq:V-generator} is bounded in absolute value by
$2\norm F_\infty(\lambda n+\gamma\binom n2)$.
It is therefore integrable in time and expectation on bounded intervals by Lemma~\ref{lem:degree-nonexplosion}, recalling that $L_t=d(V_t)$.
The same assertion applies when $F$ is bounded only on the invariant set reachable from its initial state.
More precisely, starting from a fixed $v$, every reachable coefficient function $w$ satisfies $\|w\|_\infty\leq\|v\|_\infty$, hence $|H(w,\z)|\leq\|v\|_\infty$ for $\z\in\Delta_d$.
By Lemma~\ref{lem:operator-contraction}, $\norm{V_t}_\infty\leq\norm{V_0}_\infty$ at every jump and between jumps.
Consequently $\V^{[0,1]}$ is invariant.
The construction also shows that $V$ has right-continuous paths with left limits and is stochastically continuous for the metric \eqref{eq:V-metric}.
Indeed, there are finitely many jumps on bounded intervals, and no jump occurs at a prescribed deterministic time.

\section{Proof of multivariate Bernstein duality}\label{sec:duality}

In this section we prove our main result. We first identify the two generator actions and then justify the duality identity for every continuous forward solution.
The latter step requires an integrable bound because the degree of the dual is unbounded and the coalescence rate is quadratic.

\begin{proof}[Proof of Theorem~\ref{thm:bernstein-duality}]
For fixed $v$, the function $H(v,\cdot)$ is a polynomial and belongs to the $C^2$ test class of the forward martingale problem.
For fixed $\x$, the expression $\mathcal BH(\cdot,\x)(v)$ is well defined by the two possible jump transitions in \eqref{eq:V-generator}.
For $n\geq2$, $v\in\V_n$, and $\x\in\Delta_d$, adding the identities of Propositions~\ref{prop:selection-intertwining} and \ref{prop:coalescence-intertwining} gives
\begin{equation}\label{eq:generator-duality}
\mathcal LH(v,\cdot)(\x)=\mathcal BH(\cdot,\x)(v).
\end{equation}
When $n=1$, all second derivatives of $H(v,\cdot)$ vanish and the coalescence part of $\mathcal B$ is absent, so the same identity holds.

Lemma~\ref{lem:forward-wellposed} constructs a continuous forward solution independently of duality.
Fix $v\in\V_n$ and $\x\in\Delta_d$, let $\X$ be any continuous martingale-problem solution started at $\x$, and construct $V$ independently with $V_0=v$.
We verify the hypotheses of \cite[Theorem~4.4.11 and Corollary~4.4.13]{EthierKurtz1986}, with both potential functions in that theorem equal to zero.

Write $M:=\norm v_\infty$.
For this application, $V$ may be restricted to the invariant state space $\{w\in\V:\|w\|_\infty\leq M\}$, on which $H$ is bounded.
By the contractions of Lemma~\ref{lem:operator-contraction}, every coefficient jump preserves the interval $[-M,M]$.
Nonexplosion therefore gives, on one event of probability one,
\[
\sup_{s\geq0}\sup_{\z\in\Delta_d}|H(V_s,\z)|
\leq\sup_{s\geq0}\norm{V_s}_\infty\leq M.
\]
By \eqref{eq:V-generator} and the same contractions,
\begin{equation}\label{eq:B-H-bound}
    |\mathcal BH(\cdot,\z)(V_s)|
      \leq2M\left(\lambda L_s+\gamma\binom{L_s}{2}\right),
      \qquad s\geq0,\quad\z\in\Delta_d.
\end{equation}
The generator identity gives the identical bound for $|\mathcal LH(V_s,\cdot)(\z)|$.
For $T>0$, put $\overline L_T:=\sup_{s\leq T}L_s$.
Lemma~\ref{lem:degree-nonexplosion} implies that $\Gamma_T:=M+4M\left(\lambda\overline L_T+\gamma\binom{\overline L_T}{2}\right)$
is integrable and bounds the sum of the absolute values of $H(V_r,\X_u)$ and both generator expressions, uniformly for $r\vee u\leq T$.

By localization, one can prove that for each fixed $\z$,
\[
\Big(H(V_{t},\z)-H(v,\z)-\int_0^{t}\mathcal BH(\cdot,\z)(V_s)\dd s\Big)_t
\]
is an integrable martingale.
Similarly, for each fixed $w\in\V$, since $H(w,\cdot)$ is a polynomial and the simplex is compact, 
\[
\Big(    H(w,\X_t)-H(w,\X_0)
       -\int_0^t\mathcal LH(w,\cdot)(\X_s)\dd s
\Big)_t
\]
is an integrable martingale.

All functions involved are measurable: on each fixed-degree component, the coefficient operators are linear and $H$ and its generator expressions are continuous.
The martingale identities, the integrable envelope, and \eqref{eq:generator-duality} permit application of \cite[Theorem~4.4.11  and Corollary~4.4.13]{EthierKurtz1986}, yielding \eqref{eq:semigroup-duality} for almost every $t$.
Both sides are right-continuous by the path regularity of $\X,V$ and dominated convergence, using the bound $M$.
Approaching an arbitrary time from above through the full-measure set of valid times proves the identity for every $t\geq0$.
The same proof applies when $\X_0$ has any probability law $\mu$ on $\Delta_d$ and is independent of $V$.
No uniqueness or Markov property of the arbitrary forward solution is used in this argument.
\end{proof}

Taking $H(\delta_{\kvec},\x)=\binomulti n{\kvec}\x^{\kvec}$, we obtain the following result.

\begin{corollary}[Multivariate moment duality]\label{cor:moment-duality}
Fix $t\geq0$, $n\in\N$, $\x\in\Delta_d$, and $\kvec\in\Kc_n$.
Let $\delta_{\kvec}\in\V_n$ equal one at $\kvec$ and zero elsewhere.
Then
\begin{equation}\label{eq:moment-duality}
    \E_{\x}\left[\prod_{j=1}^dX_j(t)^{k_j}\right]
    =\binomulti n{\kvec}^{-1}
    \E_{\delta_{\kvec}}[H(V_t,\x)].
\end{equation}
\end{corollary}
Expanding $H(V_t,\x)$ at its current degree $L_t$ gives \eqref{eq:moment-dualityIntro}.
Thus a forward mixed moment is computed by averaging a polynomial generated by the ancestry.

\begin{corollary}[Well-posedness of the forward martingale problem]\label{cor:forward-wellposed}
The continuous martingale problem for $(\mathcal L,C^2(\Delta_d))$ is well posed for every initial distribution on $\Delta_d$.
Its solutions have the law of the continuous Feller strong Markov process constructed in Lemma~\ref{lem:forward-wellposed}.
\end{corollary}

\begin{proof}
For deterministic initial states, existence follows from Lemma~\ref{lem:forward-wellposed}.
Its Feller family is measurable in the initial state, so mixing the initial state with any probability measure $\mu$ gives existence for that law.
Let $\X$ and $\widetilde{\X}$ be two continuous solutions with the same initial distribution $\mu$ on $\Delta_d$.
By the proof of Theorem~\ref{thm:bernstein-duality} and bounded Fubini, for every $n\geq1$, $v\in\V_n$, and $t\geq0$,
\[
\E[H(v,\X_t)]
=\int_{\Delta_d}\E_v[H(V_t,\x)]\,\mu(\dd\x)
=\E[H(v,\widetilde{\X}_t)].
\]
By Lemma~\ref{lem:bernstein-basis}, the functions $H(v,\cdot)$, with $v\in\V$, span the restrictions of all polynomials to $\Delta_d$.
Since $\Delta_d$ is compact, these polynomials are dense in $C(\Delta_d)$.
Consequently, $\X_t$ and $\widetilde{\X}_t$ have the same distribution for every $t\geq0$.
By \cite[Theorem~4.4.2 and Corollary~4.4.3]{EthierKurtz1986}, uniqueness of the one-dimensional distributions for every initial distribution implies uniqueness of the law on $C([0,\infty),\Delta_d)$.
Therefore, the continuous martingale problem is well posed.
\end{proof}

\appendix
\section{The forward martingale problem}\label{app:forward}

Now we construct a continuous Feller strong Markov solution independently of duality.
Uniqueness is deduced from duality in Corollary~\ref{cor:forward-wellposed}.

\begin{lemma}[Equivalent representation and covariance]\label{lem:forward-edge-form}
Let $D_{ij}:=\partial_i-\partial_j$.
The operator \eqref{eq:forward-generator} has the equivalent form
\begin{equation}\label{eq:forward-edge-generator}
    \mathcal Lf(\x)
    =\lambda\sum_{i<j}x_ix_jA_{ij}D_{ij}f(\x)
     +\frac\gamma2\sum_{i<j}x_ix_jD_{ij}^2f(\x).
\end{equation}
If two $C^2$ functions defined on neighbourhoods of $\Delta_d$ agree on $\Delta_d$, this formula gives the same value for both at every point of $\Delta_d$.
Moreover, let $\sigma(\x)$ be the $d\times\binom d2$ matrix whose columns are indexed by unordered pairs $\{i,j\}$, written with $i<j$, and given by
\[
\sigma^{ij}(\x)=\sqrt{\gamma x_ix_j}(\e_i-\e_j).
\]
Equivalently, its scalar entry in row $r\in[d]$ and column $(i,j)$ is
$\sigma_{r,(i,j)}(\x)=\sqrt{\gamma x_ix_j}(\delta_{ri}-\delta_{rj})$.
Then
\begin{equation}\label{eq:forward-covariance-factorization}
    \sigma(\x)\sigma(\x)^{\mathsf T}
       =a(\x):=\gamma(\diag(\x)-\x\x^{\mathsf T}),
\end{equation}which is positive semidefinite.
\end{lemma}

\begin{proof}
Using that $\sum_{i>j}x_ix_jA_{ji}\partial_{j}f=\sum_{i<j}x_ix_jA_{ij}\partial_{i}f$, $A_{jj}=0$ and $A_{ji}=-A_{ij}$, we obtain
\[
\begin{aligned}
\sum_{j=1}^d x_j(A\x)_j\partial_j f
=\sum_{i<j}x_ix_j
   \bigl(A_{ij}\partial_i f+A_{ji}\partial_j f\bigr)=\sum_{i<j}x_ix_jA_{ij}D_{ij}f.
\end{aligned}
\]
For the second-order part, expand the square of the directional derivative, use $\sum x_i=1$ and cancel a term in $\sum_{i<j}x_ix_j\partial_{ii}f=\sum_{i=1}^dx_i\partial_{ii}f\big(1-\sum_{j=1}^ix_j \big)$, to obtain
\[
\begin{aligned}
\sum_{i<j}x_ix_jD_{ij}^2f
&=\sum_{i<j}x_ix_j
  (\partial_{ii}f+\partial_{jj}f-2\partial_{ij}f)\\
&=\sum_{j=1}^d x_j(1-x_j)\partial_{jj}f
   -2\sum_{i<j}x_ix_j\partial_{ij}f\\
&=\sum_{i,j=1}^d x_i(\delta_{ij}-x_j)\partial_{ij}f.
\end{aligned}
\]
Here the two mixed partial derivatives agree because $f$ is $C^2$.
This proves \eqref{eq:forward-edge-generator}.

Now let $f_1,f_2$ be two $C^2$ extensions with $f_1=f_2$ on $\Delta_d$.
Fix $\x\in\Delta_d$ and a pair with $x_ix_j>0$.
For $|r|<\min(x_i,x_j)$, the point $\x+r(\e_i-\e_j)$ remains in $\Delta_d$.
The functions $r\mapsto f_1(\x+r(\e_i-\e_j))$ and $r\mapsto f_2(\x+r(\e_i-\e_j))$ therefore agree near zero, so their first and second derivatives at zero agree.
Consequently $D_{ij}f_1(\x)=D_{ij}f_2(\x)$ and $D_{ij}^2f_1(\x)=D_{ij}^2f_2(\x)$.
Pairs with $x_ix_j=0$ contribute zero to the formula.
This proves $\mathcal Lf_1(\x)=\mathcal Lf_2(\x)$ even on a boundary face, where some coordinates vanish.

Finally, the diagonal entries of $\sum_{i<j}\sigma^{ij}(\sigma^{ij})^{\mathsf T}$ are $\gamma x_i(1-x_i)$, and its $(i,j)$ entry for $i\neq j$ is $-\gamma x_ix_j$.
These are precisely the entries of $a$.
For every $\z\in\R^d$, the factorization gives $
\z^{\mathsf T}a(\x)\z
=\gamma\sum_{i<j}x_ix_j(z_i-z_j)^2\geq0$.
\end{proof}

\begin{remark}
This lemma separates the contributions associated with each pair of types.
The direction $\e_i-\e_j$ increases coordinate $i$ and decreases coordinate $j$ by the same scalar amount.
In the finite frequency chain \eqref{eq:intro-finite-rates}, a replacement changes these coordinates by $+1/N$ and $-1/N$.
In the diffusion representation below, this contribution has drift $\lambda X_iX_jA_{ij}(\e_i-\e_j)\dd t$ and noise $\sqrt{\gamma X_iX_j}(\e_i-\e_j)\dd B^{ij}$.
\end{remark}

The next result constructs a continuous forward solution, represents its martingale part by Brownian noise, and proves that absent types cannot reappear.
It also proves that each entry $X_j(t)$ can be seen as the proportion of type $j$ at time $t$ in the population. 

\begin{lemma}[Existence, SDE representation, and invariant faces]\label{lem:forward-wellposed}
For every $\x\in\Delta_d$, the continuous martingale problem for \eqref{eq:forward-generator} has a solution.
There is a continuous Feller strong Markov family of such solutions.
Each constructed solution admits, on an enlarged probability space if necessary, the weak representation
\begin{equation}\label{eq:forward-sde}
    \X_t=\x+\int_0^t\bvec(\X_s)\dd s
       +\sqrt\gamma\sum_{i<j}\int_0^t
          \sqrt{X_i(s)X_j(s)}(\e_i-\e_j)\dd B_s^{ij},
\end{equation}
where the $B^{ij}$ are independent standard Brownian motions.
Every face of $\Delta_d$ is invariant, including after it is first reached.
All vertices are absorbing, and $\sum_{j=1}^dX_j(t)=1$ for all $t\geq0$ almost surely.
\end{lemma}

\begin{proof}
\emph{Existence and continuity.}
First suppose $\gamma=0$.
The polynomial vector field $\bvec$ is locally Lipschitz, so its ODE has a unique solution on a maximal interval from every $\x\in\Delta_d$.
On that interval, each coordinate satisfies
\[
x_j(t)=x_j(0)\exp\left(\lambda\int_0^t(A\x(s))_j\,\mathrm ds\right),
\]
and is therefore nonnegative.
Skew-symmetry gives $
\frac{\mathrm d}{\mathrm dt}\sum_jx_j(t)
=\lambda\x(t)^{\mathsf T}A\x(t)=0$.
The solution remains in the compact simplex, which prevents finite-time blow-up.
Continuous dependence on the initial state gives a deterministic Feller strong Markov family.
This establishes existence when $\gamma=0$, including the constant process if also $\lambda=0$.

Now suppose $\gamma>0$.
The simplex has $d$ frequencies but only $d-1$ independent coordinates.
To apply \cite[Theorem~8.2.8]{EthierKurtz1986}, use the affine bijection
\[
\begin{aligned}
K_{d-1}&:=\left\{\vec{\bm y}\in[0,1]^{d-1}:\sum_{j<d}y_j\leq1\right\},\\
\theta(\vec{\bm y})&:=\left(y_1,\ldots,y_{d-1},1-\sum_{j<d}y_j\right)\in\Delta_d.
\end{aligned}
\]
For $g\in C^2(K_{d-1})$, put $\widetilde g(\x):=g(x_1,\ldots,x_{d-1})$ and define the normalized operator by
\[
\begin{aligned}
\mathcal Gg(\vec{\bm y})
&:=\gamma^{-1}(\mathcal L\widetilde g)(\theta(\vec{\bm y}))\\
&=\sum_{j<d}\beta_j(\vec{\bm y})\partial_jg(\vec{\bm y})
 +\frac12\sum_{i,j<d}y_i(\delta_{ij}-y_j)\partial_{ij}g(\vec{\bm y}),
\end{aligned}
\]
where $\beta_j(\vec{\bm y}):=(\lambda/\gamma)y_j(A\theta(\vec{\bm y}))_j$.
The factor $\gamma^{-1}$ puts the covariance into the normalization of the cited theorem.
The drift is polynomial, hence Lipschitz on $K_{d-1}$, and $\beta_j=0$ on $\{y_j=0\}$.
On $\{\sum_{j<d}y_j=1\}$, write $\x:=\theta(\vec{\bm y})$; then $x_d=0$, so
\[
\sum_{j<d}\beta_j(\vec{\bm y})
=\frac\lambda\gamma\sum_{j<d}x_j(A\x)_j
=\frac\lambda\gamma\bigl(\x^{\mathsf T}A\x-x_d(A\x)_d\bigr)=0.
\]
These are the boundary inequalities required in Theorem~8.2.8, with equality here.
That theorem gives a Feller generator equal to the closure of $\mathcal G$ on $C^2(K_{d-1})$, with polynomials as a core.
If $Q_t$ is its semigroup, then
$P_tf(\x):=Q_{\gamma t}(f\circ\theta)(\theta^{-1}(\x))$
defines a Feller semigroup on $C(\Delta_d)$ whose generator $\mathcal A$ is the closure of $\mathcal L$ on $C^2(\Delta_d)$.
Lemma~\ref{lem:forward-edge-form} ensures that another ambient extension of the same test function gives the same operator on the simplex.

By \cite[Theorem~4.2.7]{EthierKurtz1986}, the semigroup has a strong Markov realization with right-continuous paths and left limits.
One can prove its paths are continuous by Proposition~4.2.9 and Remark~4.2.10 of \cite{EthierKurtz1986}.
Indeed, one can fix $\x$ and $\varepsilon>0$ and construct a smooth function $f$ with $0\leq f\leq1$, equal to one near $\x$ and zero outside the ball of radius $\varepsilon$ about $\x$.

For $f\in C^2(\Delta_d)$, the semigroup property and the definition of the generator give
$\frac{\mathrm d}{\mathrm dt}P_tf=P_t\mathcal Af=P_t\mathcal Lf$
in the uniform norm.
Integrating yields $P_tf-f=\int_0^tP_r\mathcal Lf\,\mathrm dr$.
Let $(\mathcal F_s)$ be the natural filtration, and let $M_t^f$ denote \eqref{eq:forward-martingale-definition}.
For $0\leq s\leq t$, the Markov property and conditional Fubini give
\[
\E[M_t^f-M_s^f\mid\mathcal F_s]
=P_{t-s}f(\X_s)-f(\X_s)
-\int_s^tP_{r-s}\mathcal Lf(\X_s)\,\mathrm dr=0.
\]
Both $f$ and $\mathcal Lf$ are bounded on the compact simplex, so these are integrable martingales.
This proves existence of the required continuous martingale-problem solution.

\emph{The weak stochastic differential equation.}
Coordinate test functions $f(\x)=x_j$ give continuous martingales
\[
    \mathfrak M_j(t):=X_j(t)-x_j-\int_0^t b_j(\X_s)\dd s.
\]
They are square-integrable on bounded intervals, since $\X$ and $\bvec$ are bounded.
For $i,j\in[d]$, the quadratic test function $f(\x)=x_jx_i$ satisfies
\[
\mathcal L(x_ix_j)=b_i(\x)x_j+b_j(\x)x_i+a_{ij}(\x).
\]
Define the martingale
\[
N_{ij}(t):=X_i(t)X_j(t)-x_ix_j
-\int_0^t[b_i(\X_s)X_j(s)+b_j(\X_s)X_i(s)+a_{ij}(\X_s)]\,\mathrm ds.
\]
The product rule for the coordinate semimartingales gives
\[
\begin{aligned}
N_{ij}(t)
={}&\int_0^tX_i(s)\,\mathrm d\mathfrak M_j(s)
+\int_0^tX_j(s)\,\mathrm d\mathfrak M_i(s)+\langle\mathfrak M_i,\mathfrak M_j\rangle_t
-\int_0^ta_{ij}(\X_s)\,\mathrm ds.
\end{aligned}
\]
The stochastic integrals are continuous local martingales: their integrands are bounded and predictable because $\X$ is continuous and adapted.
Subtracting them from $N_{ij}$ shows that the last two terms form a continuous local martingale.
They also have finite variation and vanish at zero, so their difference is identically zero.
Therefore
\[
\langle\mathfrak M_i,\mathfrak M_j\rangle_t
=\int_0^ta_{ij}(\X_s)\,\mathrm ds.
\]

Together with \eqref{eq:forward-covariance-factorization}, this identity is the covariance condition for the Brownian representation theorem \cite[Theorem~5.3.3]{EthierKurtz1986}.
To use its square-matrix formulation, put $q:=\binom d2$ and $D:=\max\{d,q\}$.
Embed $\X$ and $\bvec$ in $\R^D$ by adjoining zero coordinates, and extend the $d\times q$ matrix $\sigma$ to a $D\times D$ matrix by zero rows and columns.
Extend the coefficients by zero off the embedded simplex; they are bounded Borel functions.
Every smooth test function on $\R^D$ restricts to a member of $C^2(\Delta_d)$, so the embedded process solves the corresponding continuous martingale problem.
The cited theorem supplies a $D$-dimensional Brownian motion on an extension of the probability space.
Write $\mathfrak M=(\mathfrak M_1,\ldots,\mathfrak M_d)^{\mathsf T}$,
and let $B$ be the $q$-dimensional Brownian motion whose
coordinates are the $B^{ij}$, ordered as the columns of $\sigma$.
Keeping the first $d$ state coordinates and the Brownian coordinates corresponding to the original \(q\) columns of \(\sigma\) gives
\[
\mathfrak M(t)=\int_0^t\sigma(\X_s)\,\mathrm dB_s,
\]
which, inserted into the definition of $\mathfrak M$, proves \eqref{eq:forward-sde}.

\emph{Conservation of total frequency.}
For every $\x\in\Delta_d$,
\[
\sum_{j=1}^d b_j(\x)=\lambda\x^{\mathsf T}A\x=0,
\qquad
\sum_{r=1}^d(\e_i-\e_j)_r=0\quad(i<j).
\]
Summing the coordinates of \eqref{eq:forward-sde} gives $\sum_{j=1}^dX_j(t)=\sum_{j=1}^dx_j=1$ for all $t\geq0$ almost surely.
Nonnegativity is supplied by the simplex-valued construction.

\emph{Invariant faces.}
Since $|(A\x)_j|\leq1$, we have $b_j(\x)\leq\lambda x_j$.
Consequently,
\[
    e^{-\lambda t}X_j(t)
       =x_j+\int_0^t e^{-\lambda s}\dd\mathfrak M_j(s)
        +\int_0^t e^{-\lambda s}(b_j(\X_s)-\lambda X_j(s))\dd s
\]
is a nonnegative supermartingale.
The stochastic integral is a true martingale on every bounded interval because $a_{jj}(\x)=\gamma x_j(1-x_j)\leq\gamma/4$.
If $x_j=0$, then $\E_{\x}[e^{-\lambda t}X_j(t)]\leq0$, so $X_j(t)=0$ almost surely for every fixed $t$.
Taking all rational times and using continuity gives this assertion simultaneously for all $t\geq0$.
Now let $\tau_j:=\inf\{t\geq0:X_j(t)=0\}$.
On $\{\tau_j<\infty\}$, the strong Markov property applies to the process started at $\X_{\tau_j}$, whose $j$th coordinate is zero.
For a possibly infinite $\tau_j$, apply this property at $\tau_j\wedge m$ on $\{\tau_j\leq m\}$ and then let $m\to\infty$.
Thus $X_j(\tau_j+t)=0$ for every $t\geq0$, almost surely on $\{\tau_j<\infty\}$.
Applying this to each absent coordinate proves invariance of every face and absorption at the vertices.
\end{proof}

\section*{Acknowledgments}

I am very grateful to Laura Eslava (Instituto de Investigaciones en Matem\'aticas Aplicadas y en Sistemas, Universidad Nacional Aut\'onoma de M\'exico) and Mar\'ia Clara Fittipaldi (Facultad de Ciencias, Universidad Nacional Aut\'onoma de M\'exico) for their insightful discussions, which were instrumental in establishing the duality result in a simpler case.

Most of this project was done during my postdoctoral positions at University Duisburg-Essen, Germany and Centro de Investigaci\'on en Matem\'aticas, M\'exico. 
I acknowledge support from the Deutsche Forschungsgemeinschaft through DFG-SPP-2265.

\providecommand{\bysame}{\leavevmode\hbox to3em{\hrulefill}\thinspace}

\end{document}